\documentclass [12pt]{amsart}
\newif\ifpdf
\ifx\pdfoutput\undefined
	\pdffalse
\else
	\pdfoutput=1
	\pdftrue
\fi
\usepackage{enumerate}
\usepackage{amssymb}
\ifpdf
	\usepackage[pdftex]{graphicx}
	\DeclareGraphicsExtensions{.pdf}
	\usepackage{hyperref}
\else
	\usepackage{graphicx}
	\DeclareGraphicsExtensions{.eps}
	\usepackage{hyperref}
\fi

\newcommand{\ZZ}{\mathbb Z}
\newcommand{\PP}{\mathbb P}

\newcommand{\CC}{\mathbb C}

\newcommand{\bbA}{\mathbb A}
\newcommand{\BB}{\mathbb B}

\newcommand{\GL}{\mathop {\rm GL}\nolimits}

\newcommand{\Pic}{\mathop {\rm Pic}\nolimits}
\newcommand{\Supp}{\mathop {\rm Supp}\nolimits}

\newtheorem{thm}{Theorem}[section]
\newtheorem{thm0}{Theorem}
\newtheorem{cor}[thm]{Corollary}
\newtheorem{cor0}[thm0]{Corollary}
\newtheorem{prop}[thm]{Proposition}
\newtheorem{lem}[thm]{Lemma}

\newtheorem*{claim}{Claim}
\theoremstyle{definition}

\newtheorem*{ack}{Acknowledgment}
\theoremstyle{remark}
\newtheorem{rem}[thm]{Remark}
\newtheorem{rem0}[thm0]{Remark}

\renewcommand{\thesubparagraph}{\theparagraph.\@arabic\c@subparagraph}
\title{Nodal degenerations of plane curves and Galois covers}

\author[E. Artal]{Enrique Artal Bartolo}

\author[J.I. Cogolludo]{Jos\'e Ignacio Cogolludo}
\address{Departamento de Matem\'aticas\\
Universidad de Zaragoza\\
Campus Plaza San Francisco s/n\\
E-50009 Zaragoza SPAIN}
\email{artal@unizar.es,jicogo@unizar.es}

\author[H. Tokunaga]{Hiro-o Tokunaga}
\address{Department of Mathematics \\
Tokyo Metropolitan University \\
Minamiohsawa Hachoji  \\
192-0357  Tokyo JAPAN}
\email{tokunaga@comp.metro-u.ac.jp}
\thanks{First and second authors are partially supported by
BFM2001-1488-C02-02.}
\begin{document}
\begin{abstract}
Globally irreducible nodes 
(i.e. nodes whose branches belong to the same irreducible component) 
have mild effects on the most common topological invariants of an algebraic curve.
In other words, adding a globally irreducible node (simple nodal degeneration)
to a curve should not change them a lot. In this paper we 
study the effect of nodal degeneration of curves on fundamental groups 
and show examples where simple nodal degenerations produce non-isomorphic 
fundamental groups and this can be detected in an algebraic way by means of Galois coverings.
\end{abstract}
\maketitle
\section*{Introduction}

Let $\{C_t\}_{t\in \Delta}$, $\Delta := \{ t \in \CC \mid |t| < 1\}$ be a family of curves such that $\{C_t\}_{t\in \Delta^*}$ is equisingular and $C_0$ is reduced. Under these
conditions, there is a surjection 
$h_t : \pi_1(\PP^2\setminus C_0) \to \pi_1(\PP^2\setminus C_t)$, see \cite[p.~121]{dimca}.
The interesting case arises when the family is not equisingular at $C_0$; the complexity of the family may be measured by the \emph{new} singularities of $C_0$. The simplest case appears when $C_0$ has the same number of irreducible components as $C_t$, $t\neq 0$, and the new
singular points are nodal points; we call them \emph{nodal degenerations}. A nodal 
degeneration is called \emph{simple} if only one node has been added.

For example, let $C_0$ be a nodal irreducible curve of degree $d$; it is easy to construct a
family where  $C_t$, $t\neq 0$, is smooth. In this case, $h_t$ is an isomorphism. This approach appears
already in \cite{zariski1}; the last statement was known as the Zariski problem. Zariski \cite{zariski1} showed
that $h_t$ is an isomorphism using a result of Severi \cite{sev:21} about the irreducibility of the moduli space of
nodal curves. Severi's proof was not right and a correct one was given by Harris~\cite{har:84}. Meanwhile,
the works of Fulton~\cite{ful:80} and Deligne~\cite{del:80} provide a proof of Zariski problem independent
of Severi's result. Deligne's approach was rather topological and one could think that it could be extended to
other nodal degenerations. Nevertheless, classical results by Zariski provide examples where $h_t$ is not an isomorphism. Using generic plane sections of discriminantal complements, Zariski~\cite{zariski2} showed that the fundamental
group of a sextic with six cusps and four nodes is isomorphic to the braid group of four strings on the $2$-sphere; there is a nodal deformation where the generic element is a 
sextic with six cusps and the fundamental group is $\ZZ/2*\ZZ/3$. In \cite{tokunaga4}, 
this curve provides an example of a simple nodal deformation with non-isomorphic $h_t$ 
(the curve $C_t$ has six cusps and three nodes); the fact that $h_t$ is not an isomorphism is shown using
finite Galois covers. Other examples were also given in \cite{tokunaga4}.

It was already shown by Zariski~\cite{zr:31}
that for nodal-cuspidal curves, nodes have no effect on the computation of the Alexander polynomial of a curve. The general case was shown by Esnault~\cite{es:82} and Libgober~\cite{li:82}. The generalization of these arguments by Libgober~\cite{li:01} showed that it is also the case for non-coordinate characteristic varieties. Experience also shows
that in most cases $h_t$ is an isomorphism for (simple) nodal deformations.

In this article we study both topologically and algebraically how to find new 
(and essentially different) families where
nodal degenerations have an effect on the fundamental group. We study certain families of curves such that the kernels of $h_t$ are 
non-trivial by using Galois covers having dihedral groups $D_{2n}$ and a semi-direct product
of $\ZZ/k\ZZ$ and $(\ZZ/2\ZZ)^{\oplus (k-1)}$, denoted by $G(k)$, as their Galois group,
introduced in \S\ref{gal-cover}. The \emph{generic} cases of these families are nodal
degenerations and fundamental groups can be computed. We also give a sufficient condition 
for a degeneration to satisfy the isomorphism statement in terms of cuspidal degenerations, see Remark~\ref{rem-bifam}.
One striking feature is the influence of the rationality of the irreducible components in the existence of Galois coverings, see Theorems~\ref{thm:main1} and~\ref{thm:main2}. This
condition appears in a natural way when studying Galois covering and it is confirmed by
the actual computations of fundamental groups. This rationality property also appears in classical examples, though Lemma~\ref{lem:preimD}\eqref{lem:preimDii} and Proposition~\ref{prop-grupo2}\eqref{prop-grupo2c}
seem to suggest that it is not essential for the non-isomorphism property of nodal
degenerations.

\medbreak
Let us describe the curves we consider in this article:

\begin{description}
\item[Type I] Let $D$ be an irreducible plane curve of degree $d$ such that
$D$ has an ordinary $(d-2)$-ple point at $P$.
Let $L_1$ and $L_2$ be lines through $P$ such that either $L_i$ is tangent to a smooth
point $P_i\in D$ or $L_i$ passes through a double point $P_i\neq P$ of type $\bbA_{2r}$.
Let us denote $C := L_1 + L_2 + D$.

\medbreak\item[Type II] Let $D$ be an irreducible plane curve of degree $d$ such that 
$D$ has an ordinary $(d-k)$-point at $P$. Let $L_1$ and $L_2$ be lines through $P$ such 
that both of $L_i$ $(i = 1, 2)$ are tangent to $D$ at inflection points $P_i$ with maximal 
order of contact (i.e. $k$). Let us denote $C := L_1 + L_2 + D$.
\end{description}

\medbreak
It is easy to prove the existence of curves of Type I and II, since we can fix
the points $P,P_1,P_2$ in a coordinate system and then impose the necessary linear 
equations required by both types.

\medbreak

In this article, we prove the following two statements (as for $D_{2n}$ and 
$G(k)$-covers see \S\ref{gal-cover}).

\medbreak

\begin{thm0}\label{thm:main1}
Let $C$ be a curve of Type I.
\begin{enumerate}[\rm(i)]
\item If $D$ is rational (i.e., the normalization of $C$ is isomorphic to $\PP^1$), 
then there exists a $D_{2n}$-cover $\pi: S \to \PP^2$ branched $2(L_1+L_2) + nD$ 
for any $n \ge 3$.

\item If $D$ is not rational, then there exists no $D_{2n}$-cover $\pi: S \to \PP^2$ 
branched at $2(L_1+ L_2) + nD$ for any odd $n \ge 3$.
\end{enumerate}
\end{thm0}

\begin{thm0}\label{thm:main2}
Let $C$ be a curve of Type II.
\begin{enumerate}[\rm(i)] 
\item\label{thm:main2i} If $D$ is rational
then there exists a $G(k)$-cover
$\pi: S \to \PP^2$  such that it is branched at
$k(L_1+L_2) + 2D$ and $D(S/\PP^2)$ is a cyclic cover branched at $k(L_1 + L_2)$.

\item\label{thm:main2ii} If $D$ is not rational, then there exists no $G(k)$-cover 
$\pi: S \to \PP^2$  such that it is branched at
$k(L_1+ L_2) + 2D$ and $D(S/\PP^2)$ is a cyclic cover branched at $k(L_1 + L_2)$.
\end{enumerate}
\end{thm0}

A straightforward consequence is the following:

\begin{cor0}\label{cor:main}
Let $C_t, t\in\Delta$ be a family of curves such that $C_t$, $t\neq 0$ is of type 
I or type II.
If $C_0$ is rational and $C_t$ ($t\neq 0$) is non-rational, then
the kernel of $h_t : \pi_1(\PP^2 \setminus C_0) \to \pi_1(\PP^2\setminus C_t)$ is non-trivial.
\end{cor0}

Since families of types I and II are non-empty for each degree, Corollary~\ref{cor:main}
provides an infinite family of examples of simple nodal degenerations such that the
kernel of $h_t$ is non-trivial. Using Theorem~\ref{thm:main1}
and~\ref{thm:main2}, and the results of \cite{ACT04b}, it is possible to construct rational
curves of type I and II having only nodes as singularities outside $P,P_1,P_2$.
In the case of curves of type I the existence of simple nodal degenerations follows from
the $T$-property of the equisingular families, see~\cite{shu:96}. 
In the case of curves of type II, examples of simple degenerations
can be found at least for small degrees.

\begin{rem0}\label{rem:deg}
Theorems~\ref{thm:main1} and \ref{thm:main2} are also true for the following degenerations
of Type I and II, i.e. allowing a point $P_i$ to be infinitely near to $P$. Proofs of such
results follow along the lines of those presented here.
\end{rem0}

\begin{ack}
Part of this work was done during the last author's visit to Ruhr Universit\"at Bochum under the support from Professor Huckleberry. He thanks Professor Huckleberry
for his hospitality and encouragement. Also he thanks  Universidad de Zaragoza for its nice
atmosphere.
\end{ack}


\section{Preliminaries}\label{prel}

Throughout this article, $\Sigma_n$ denotes a Hirzebruch surface of degree $n$ and
$p_n:\Sigma_n\to\PP^1$ the ruling. It is well known that
$\Pic(\Sigma_n) \cong \ZZ \Delta_{0, n} \oplus \ZZ F_n$, where $\Delta_{0, n}$ denotes the section with $\Delta_{0,n}^2 = -n$ and $F_n$ is a fiber. Let $C$ be a curve of either type 
I or II described in the introduction.
Let $q_P : \Sigma_1 \to \PP^2$ be the blowing-up at $P$. We denote the proper transforms 
of $D$, $L_1$ and $L_2$ by ${\overline D}$, ${\overline L}_1$, and ${\overline L}_2$,
respectively. Note that $\Delta_{0,1}$ is the exceptional divisor of $q_P$. 
Let $f_k : \Sigma_k \to \Sigma_1$
be a $k$-cyclic cover branched at $k({\overline L}_1 + {\overline L}_2)$,
$k=2$ if $D$ is of type I. We denote by $\sigma_k$ a generator of the cyclic group
of covering transformations. Note that such a $k$-cyclic cover is unique, as
$\pi_1(\Sigma_1\setminus (\overline{L}_1 
\cup \overline{L}_2)) \cong  \pi_1(\PP^1\times \CC^{*}) \cong \ZZ$.

\begin{lem}\label{lem:preimD} 
Let $G_1,\dots,G_\lambda$ be the irreducible components of $f_k^*{\overline D}$.

\begin{enumerate}[\rm(i)]
\item If $D$ is rational, then $\lambda=k$ and $G_i$ is linear equivalent to 
$\Delta_{0, k} + d F_k$.

\item\label{lem:preimDii} If $D$ is not rational then $\lambda | k$ and $\lambda<k$.
\end{enumerate}
\end{lem}

\begin{proof}  Since $D$ has an ordinary $(d-k)$-ple point at $P$,  
then
$\overline D \sim k\Delta_{0,1} + dF_1$.

\begin{enumerate}[(i)]
\item Suppose that $D$ is rational.  By the assumption of the intersection of $D$ and 
$L_1\cup L_2$ the covering $f_{k|}:f_k^*{\overline D}\to\overline D$ is unramified and 
then $f_k^*{\overline D}$ consists of $k$ irreducible components. 

Let us prove it; let $\widetilde D$ and $\widetilde{\mathcal D}$ be the normalizations of
$\overline D$ and $f_k^*{\overline D}$, respectively. The restriction of $f_k$ induces
a covering $\tilde f_k:\widetilde{\mathcal D}\to\widetilde D$ between compact smooth curves. 
This covering is unramified outside the preimage on $\widetilde D$ of
$\overline D\cap(\overline L_1\cup\overline  L_2)$. Since $D\cap L_i=\{P,P_i\}$ and
$D$ is locally irreducible at $P_i$ (which is why the index of $P_i$ is even in type I curves),
such preimage consists only of two points 
$\widetilde P_1$ and $\widetilde P_2$ (the preimages of $P_1$ and $P_2$ resp.). 
Thus the covering $\tilde f_k$ is controlled by the monodromy 
$\tau:\pi_1(\widetilde D\setminus\{\widetilde P_1,\widetilde P_2\})=\pi_1(\CC^*)=\gamma\ZZ \to\ZZ/k\ZZ$, where $\gamma$ is a meridian in $\widetilde D$ around say $\widetilde P_1$.
Since $k=(\overline D\cdot\overline L_i)_{P_i}$, one has that $\tau(\gamma)=0$, thus 
$\tau$ is trivial and hence $\widetilde{\mathcal D}$ consists of $k$ connected components.

\medbreak
Since $G_i^{\sigma_k} = G_{i+1}$, $F_k^{\sigma_k} = F_k$, and 
$\Delta_{0, k}^{\sigma_k} = \Delta_{0, k}$ we have that: 
$$
G_i\cdot F_k=G_j\cdot F_k
\text{ and }
G_i\cdot \Delta_{0, k}=G_j\cdot \Delta_{0, k}
$$
for any $i,j\in \{1,...,k\}$.
Thus
$$
k = f_k^*{\overline D}\cdot F_k = \left(\sum_{i=1}^k G_i\right) \cdot F_k=
k\,G_i\cdot F_k,
$$ 
i.e., $G_i\cdot F_k = 1$ ($i = 1, \dots, k$). Also,
$$
k(d-k)=f_k^*{\overline D} \cdot \Delta_{0,k}= \left(\sum_{i=1}^k G_i\right) \cdot \Delta_{0, k} =
k\,G_i\cdot \Delta_{0, k},
$$
and hence $\Delta_{0, k}\cdot G_i  = d - k$ ($i = 1, \dots, k$). 
Therefore $G_i \sim \Delta_{0, k} + dF_k$.

\medbreak\item 
Since $\{ G_1,\dots, G_\lambda\}$ is a $\ZZ/k\ZZ$-orbit, one always has that $\lambda | k$. 
Therefore, in order to finish the proof it is enough to show that $\lambda\neq k$.
Suppose that $f^*{\overline D}$ consists of $k$ irreducible components. 
As we have seen, this implies that $G_i\cdot F_k = 1$ ($i = 1, \dots, k$), that is 
each irreducible component $G_i$ is a section and hence it is rational, but this would
force $\tilde D$ (and hence $D$) to be rational which contradicts the hypothesis. 
\end{enumerate}
 \end{proof}

\section{A simple example of nodal degeneration}
\label{sec-example}
\begin{figure}[ht]
\begin{center}
\includegraphics[width=5cm,bb=0 0 561 340]{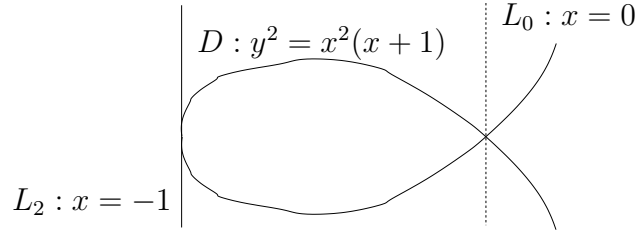}
\end{center}
\caption{Affine real picture of $C$}
\label{fig-af}
\begin{picture}(0,0)
\put(-65,110){$D:y^2=x^2 (x+1)$}
\put(50,120){$L_0:x=0$}
\put(-135,50){$L_2:x=-1$}
\end{picture}
\end{figure}

We consider an example which is a \emph{degenerated} case
of curves of type I as in Remark~\ref{rem:deg}. Let $D$ be an irreducible nodal cubic;
let $L_1$ be the tangent at an inflection point $P$
and let $L_2$ be the unique ordinary tangent to $D$ passing
through $P$. Let $C=D\cup L_1\cup L_2$. If we identify $L_1$ as
the line at infinity, where $P$ is the point at infinity
of vertical lines, $C$ has a nice real picture, see Figure~\ref{fig-af}
with affine equations.

Let $L_0$ be the line joining $P$ and the node.
The fundamental group $\pi_1(\PP^2\setminus(C\cup L_0))$
can be computed using the fibered version of Zariski-van Kampen
Theorem, see~\cite[Prop.~3.4]{ACC03}. In order to compute
this group we consider the braid monodromy of the affine part of $D$
with respect to the projection $(x,y)\mapsto x$. Since $x=0,-1$
are the only non-transversal vertical lines, braid monodromy
is a mapping $\rho:\pi_1(\PP^1\setminus\{0,-1\};-\frac{1}{2})\to\BB_2$,
where $\BB_2$ is the braid group on two strings, which is isomorphic
to $\ZZ$, having as generator the standard half-twist $\sigma_1$.
Let $\alpha,\beta$ be meridians around $0,-1$, respectively,
generating $\pi_1(\PP^1\setminus\{0,-1\};-\frac{1}{2})$ as
a free group, see~\cite[Def.~1.12]{ACC03} for a definition of meridian. 
It is easily seen that $\rho(\alpha)=\sigma_1^2$
and  $\rho(\beta)=\sigma_1$.

Let $L$ be the vertical affine line $x=-\frac{1}{2}$.
The group $\pi_1(\PP^2\setminus(C\cup L_0))$ is obtained
as a quotient of the group $\pi_1(L\setminus D)$, freely generated
by two (suitably chosen) meridians $m_1,m_2$ of the points $L\cap D$ in $L$.
The relations are of type $m=m^{\rho(u)}$, where $m\in\pi_1(L\setminus D)$,
$u\in\pi_1(\PP^1\setminus\{0,-1\};-\frac{1}{2})$ and the superscript
means the action of $\BB_2$ on $\pi_1(\PP^2\setminus(C\cup L_0))$
defined by $m_1^{\sigma_1}=m_2$ and $m_2^{\sigma_2}=m_2 m_1 m_2 ^{-1}$.
Choosing meridians $a,b$ (of $L_0,L_2$, respectively)
which are suitable liftings of $\alpha,\beta$, we obtain:
\begin{equation*}
\begin{split}
\pi_1(\PP^2\setminus(C\cup L_0))=&
\langle m_1,m_2,a,b\mid\\
&m_1^a=m_2 m_1 m_2 ^{-1}, [m_2 m_1, a]=1,
m_1^b=m_2, [m_2 m_1, b]=1
\rangle.
\end{split}
\end{equation*}

Applying the techniques of the standard Zariski-van Kampen theorem
in order to obtain $\pi_1(\PP^2\setminus C)$ it is enough to \emph{kill}
the meridian $a$ of $L_0$:
\begin{equation*}
\begin{split}
\pi_1(\PP^2\setminus C)=&
\langle m_1,m_2,b\mid [m_2, m_1]=1,
m_1^b=m_2, [m_2 m_1, b]=1
\rangle
\\
=&\langle m_1,b\mid [m_1, b^2]=1,
[m_1, b m_1 b]=1
\rangle.
\end{split}
\end{equation*}
By setting $m:=m_1 b$, one obtains that:
$$
\pi_1(\PP^2\setminus C)=\langle m,b\mid [m, b^2]=1,
[m^2, b]=1\rangle.
$$
Note that $m^2,b^2$ are central elements and since
the quotient by the subgroup generated by these elements
is a free product of two copies of $\ZZ/2\ZZ$, it can be seen
that they generate the center and that $\pi_1(\PP^2\setminus C)$
is not Abelian.

\begin{figure}[ht]
	\centering
	\includegraphics[width=5cm,bb=0 0 581 366]{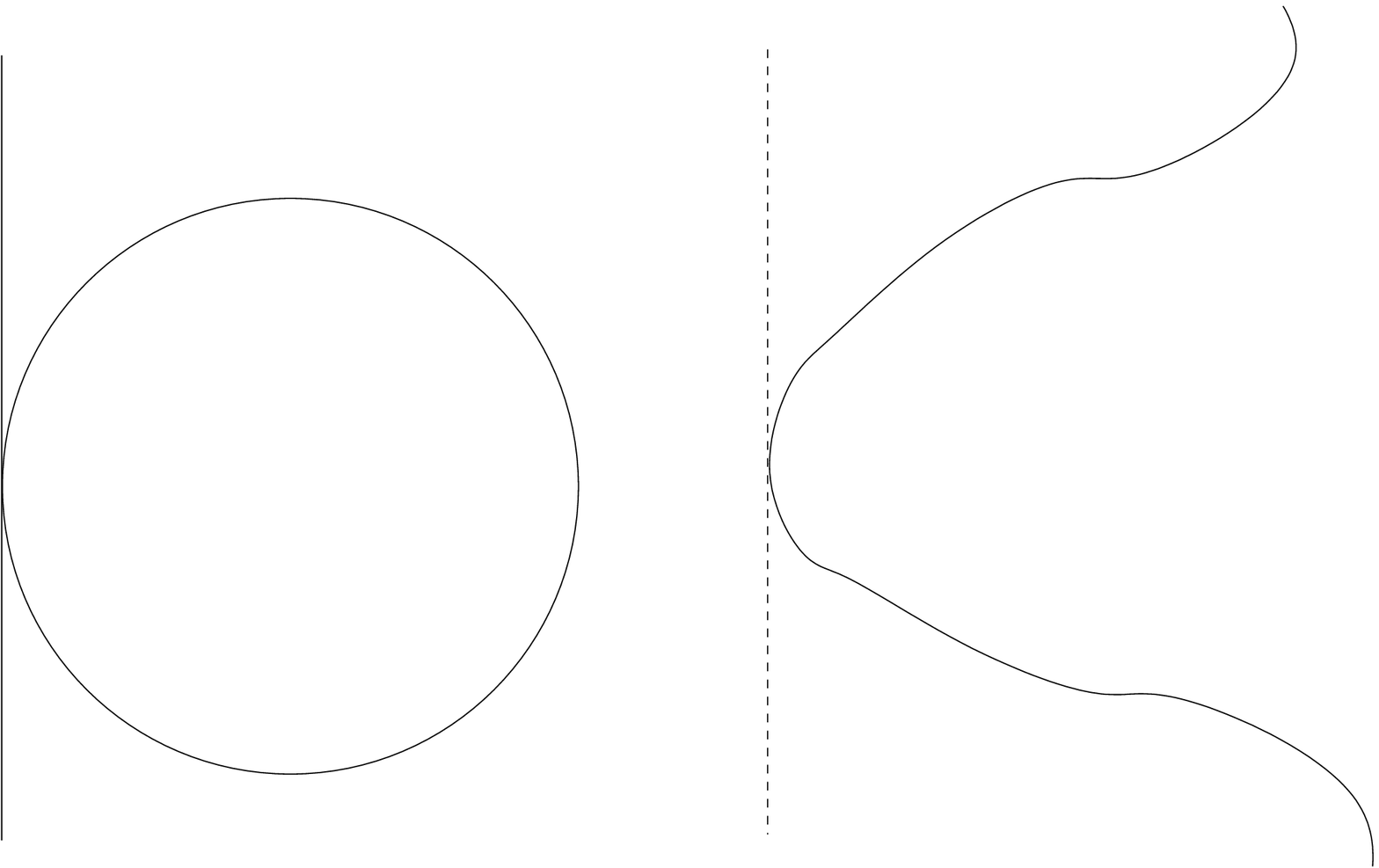}
	\caption{Real affine part of $C_t$ in $\PP^2\setminus L_1$}
	\label{fig-smooth}
\begin{picture}(0,0)
\put(-30,110){$D_t$}
\put(-7,50){$L_0$}
\put(-85,50){$L_2$}
\end{picture}
\end{figure}

We can consider this curve in a family $C_t$, $t\in\Delta$, such that
$C_t=D_t\cup L_1\cup L_2$, where $D_t$ (for
$t\neq 0$) is a smooth cubic  having an inflection
point at $P$, the line $L_1$ is tangent to $D_t$ at $P$, and $L_2$ is a tangent 
line to $D_t$ passing through $P$, see Figure~\ref{fig-smooth}.
It is easily seen that in this case $\pi_1(\PP^2\setminus C_t)$
is isomorphic to $\ZZ^2$. This is the simplest example of a simple 
nodal degeneration with non-isomorphic fundamental groups.

\begin{rem}\label{rem-bifam}
Let us suppose that we have a family $C_{t,s}$, $(t,s)\in\Delta^2$
such that the whole family is equisingular outside a small polydisk $B$
in $\PP^2$ where the curves behave as in Figure~\ref{fig-deg}.
\begin{figure}[ht]
	\centering
	\includegraphics[width=14cm,bb=0 0 627 158]{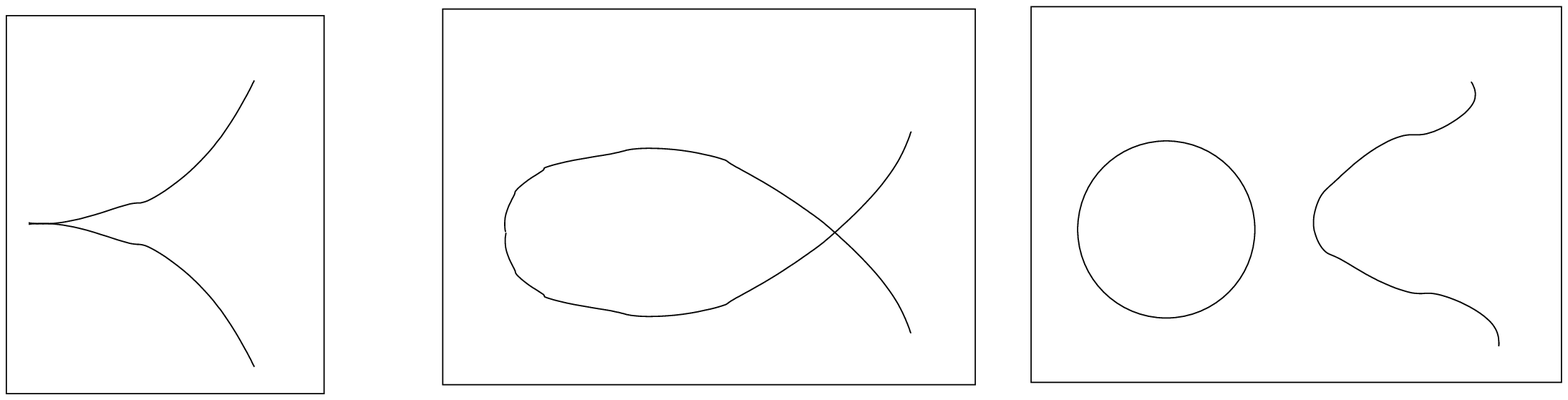}
	\caption{Bifamily}
	\label{fig-deg}
\begin{picture}(0,0)
\put(100,30){$C_{t,s}, s\neq 0$}
\put(-40,30){$C_{t,0}, t\neq 0$}
\put(-170,30){$C_{0,0}$}
\end{picture}
\end{figure}
Let us apply the classic Zariski-van Kampen theorem to this family.
The fundamental group $\pi_1(\PP^2\setminus C_{t,s})$ is the quotient
of $\pi_1(L\setminus C_{t,s})$, $L$ a \emph{generic} transversal line,
by the system of relations given by the braid monodromy action. 
We can choose $L$ intersecting $B$ and we have two special free 
generators $m_{t,s},n_{t,s}$ of $\pi_1(L\setminus C_{t,s})$.
Braid monodromy relations outside $B$ are the same ones for each $(t,s)$.
Inside $B$, we have the following relations:
\begin{itemize}
\item  $m_{t,s}=n_{t,s}$, if $s\neq 0$.
\item  $m_{t,0}=n_{t,0}$ and $[m_{t,0},n_{t,0}]=1$, if $t\neq 0$ (note that
the second relation is a consequence of the first one).
\item  $m_{0,0} n_{0,0} m_{0,0}=n_{0,0} m_{0,0} n_{0,0}$.
\end{itemize}
We deduce that outside $(0,0)$, the groups $\pi_1(\PP^2\setminus C_{t,s})$ are isomorphic
and we obtain simple nodal deformations where the surjection of fundamental groups
is bijective. Note that we can construct bifamilies in this way where
$D_{t,0}\cup L_1=D_0\cup L_1$, $t\neq 0$. However, such a bifamily is not 
possible for $C$ because of the existence of the vertical tangents $L_1$ and $L_2$.
\end{rem}



\section{Some Galois covers}\label{gal-cover}

Let $G$ be a finite group with normal subgroup $H$. Let $X$ and $Y$  be normal
projective varieties. Let $\pi : X \to Y$ be a finite surjective morphism.
Under these conditions, the rational function field $\CC(X)$ of $X$ is regarded
as a finite field extension of $\CC(Y)$. We call $X$ a \emph{$G$-cover} of $Y$ if the
field extension $\CC(X)/\CC(Y)$ is Galois and its Galois group is isomorphic to $G$.
Let $H$ be a normal subgroup of $G$, let $\pi : X \to Y$  be a $G$-cover and let 
$\CC(X)^H$ be the fixed field of $H$. We define a normal variety $D_H(X/Y)$ to be
the $\CC(X)^H$-normalization of $Y$. Note that there are canonical morphisms
$\beta_{1,H}(\pi) : D_H(X/Y) \to Y$ and  $\beta_{2,H}(\pi) : X\to D_H(X/Y)$.
The former is a $G/H$-cover of $Y$ and the latter is a $H$-cover of $D_H(X/Y)$.
The branch locus  of  $\pi: X \to Y$, denoted by $\Delta(X/Y)$ or
$\Delta_{\pi}$ is the subset of $Y$ given by
\[
\Delta_{\pi}:=\{ y\in Y \mid \mbox {\rm $\pi$ is not a local isomorphism at $y$}\}.
\]
It is well known that $\Delta_{\pi}$ is an algebraic subset of codimension $1$
if $Y$ is smooth, see~\cite{zariski4}. Suppose that $Y$ is smooth and let
$\Delta_{\pi} =  B_1+ \dots + B_r$ be the irreducible decomposition. We say
$\pi : X \to Y$ is branched at $e_1B_1+ \cdots +e_rB_r$ if
the ramification index along $B_i$ is $e_i$.

\subsection{$D_{2n}$-covers}\mbox{}

Let $D_{2n}$ be the the dihedral group of order $2n$. 
To represent $D_{2n}$, we use the notation
\[
D_{2n} = \langle \sigma, \, \tau \mid \sigma^2 = \tau^n = (\sigma \tau)^2 = 1 \rangle,
\]
and fix it throughout this article. If we put $H = \langle \tau \rangle$, then the variety
$D_{\langle \tau \rangle}(X/Y)$ introduced as above is a double cover of $Y$, and 
$X$ is an $n$-cyclic cover of $D_{\langle \tau \rangle}(X/Y)$

A concrete method to deal with $D_{2n}$-covers was developed in~\cite{tokunaga1}.   
In this note, we only need the following simple form given in~\cite{tokunaga3}. 
As for details on $D_{2n}$-covers see \cite{tokunaga1}.

\begin{prop}\label{prop:sufdi}
Let $f : Z \to Y$ be a double cover between smooth projective varieties $Y$ and $Z$. 
Let $\sigma_f$ be the covering transformation determined by $f$. Let $D$ be a reduced 
divisor on $Z$ such that 

\begin{enumerate}[\rm(i)]
\item $D$ and $\sigma_f^* D$ have no common components,
\item if we let $D = \sum_ia_iD_i$ be the irreducible decomposition, then 
$a_i > 0$ $\forall i$ and the greatest common divisor of $a_i$'s and $n$ is $1$, and
\item $D$ is linearly equivalent to $\sigma_f^*D$.
\end{enumerate}

Then for any $n \ge 3$, there exists a $D_{2n}$-cover $\pi: X \to Y$ such that 
\begin{enumerate}[\rm(a)]
\item $D_{\langle\tau\rangle}(X/Y) = Z$, 
$f = \beta_{1, \langle \tau \rangle}(\pi)$
and 
\item the branch locus of $\beta_{2, \langle\tau\rangle}(\pi)$ is contained in 
$\Supp(D + \sigma_f^*D)$, i.e. $\Delta_{\pi} \subset \Delta_f \cup f(\Supp(D))$.
\end{enumerate}
\end{prop} 


As for a necessary condition for the existence of $D_{2n}$-covers, we have the following.

\begin{prop}\label{prop:necdi} 
Let $\pi : X \to Y$ be a $D_{2n}$-cover such that both $Y$ and
$D_{\langle\tau\rangle}(X/Y)$ are smooth. Let $\sigma$ be the covering transform of
$\beta_{1,\langle\tau\rangle}(\pi)$. Then there exist three effective divisors $D_1$, 
$D_2$ and $D_3$ on $D_{\langle\tau\rangle}(X/Y)$ satisfying the following conditions:

\begin{enumerate}[\rm(i)]
\item $D_1$ and $\sigma^*D_1$ have no common components.
\item If we let $D_1 = \sum_ia_i D_{1,i}$ be the irreducible decomposition, then 
$0 < a_i \le (n - 1)/2$ for every $i$.
\item $D_1 + nD_2 \sim \sigma^*D_1 + nD_3$ and $D_2 + \sigma^* D_2 \sim D_3 + \sigma^*D_3$.
\item $\Supp(D_1 + \sigma^*D_1) = \Delta_{\beta_{2,\langle\tau\rangle}(\pi)}$.
\end{enumerate}
\end{prop}

For a proof see \cite[\S2]{tokunaga1}.

\begin{cor}\label{cor:branch} 
Let $D$ be an irreducible component of 
$\beta_{1, \langle\tau\rangle}(\pi)(\Delta_{\beta_{2, \langle\tau\rangle}(\pi)})$.
Then $\beta_{1, \langle\tau\rangle}(\pi)^*D$ is in the form $D' + \sigma^*D'$ for some irreducible divisor 
on $D_{\langle\tau\rangle}(X/Y)$. In other
words, $\beta_{1, \langle\tau\rangle}(\pi)$ is not branched along any irreducible divisor $D$ with $D = \sigma^*D$.
\end{cor}

\subsection{Certain $\ZZ/k\ZZ\ltimes(\ZZ/2\ZZ)^{k-1}$-covers}\mbox{}

Let $k$ be an integer $\ge 3$. Let $G(k)$ be a group generated by $\sigma, \tau_1,\dots, \tau_{k-1}$ with relations
\[
\sigma^k = \tau_1^2 = \cdots = \tau_{k-1}^2 = 1,  \quad \sigma\tau_i\sigma^{-1} = \tau_{i+1}
\text{ where }\tau_k:=\prod_{j=1}^{k-1}\tau_j.
\]
The group $G(k)$ is a semi-direct product of $\ZZ/k\ZZ$ and $(\ZZ/2\ZZ)^{\oplus k-1}$.
In fact, $G(k)$ is a subgroup of index $2$
of the wreath product $\ZZ/k\ZZ$ and $\ZZ/2\ZZ$) and there is a faithful representation
$\rho_k : G(k) \to \GL (k, \ZZ)$ given by
$$
\sigma \mapsto  
\begin{pmatrix}
0 & 0 & \cdots & 0 & 1 \\
1 & 0 & \cdots & 0 & 0 \\
0 & 1 &           &    & 0 \\
0 &  0 & \ddots &  & \vdots \\
0 & 0 & \cdots & 1 & 0 
\end{pmatrix},
\quad 
\tau_1 \mapsto
\begin{pmatrix}
-1 & 0 & 0 & \cdots & 0 & 0 \\
0 & -1 & 0 & \cdots & 0 & 0 \\
0 & 0 & 1 & \cdots & 0 & 0 \\
&&& \ddots && \\
0 & 0 & 0 & \cdots & 1 & 0 \\
0 & 0 & 0 & \cdots & 0 & 1 \\
\end{pmatrix},
$$
and
$\tau_{i+1}:=\sigma \tau_i \sigma^{-1}, i=1,\dots,k-1$.
One can define a right action of $G(k)$ on $\CC(x_1,\dots, x_k)$ given by
\begin{eqnarray*}
(x_1,\dots, x_k)^{\sigma} &=& (x_2, x_3, \dots, x_k, x_1) \\
(x_1,\dots, x_k)^{\tau_i} & = & 
(x_1, \dots,x_{i-1}, x_i^{-1}, x_{i+1}^{-1},x_{i+2}, \dots, x_k)
\quad (i = 1,\dots, k-1).
\end{eqnarray*}

Let us introduce the birational coordinate change
\[
y_i = \frac {1 + x_i}{1-x_i} \quad (i = 1,\dots, k).
\]
With these new coordinates, we have
\begin{eqnarray*}
(y_1,\dots, y_k)^{\sigma} &= & (y_2,y_3,\dots, y_k, y_1) \\
(y_1,\dots, y_k)^{\tau_i} &= & (y_1,\dots,y_{i-1},-y_i, -y_{i+1},y_{i+2},\dots, y_k)  
\quad (i = 1,\dots, k-1).
\end{eqnarray*}

\begin{rem}\label{rem-versal}
The action of $G(k)$ on $\CC(x_1, \dots, x_k)$ is versal, see~\cite{bannai},
and so is on $\CC(y_1,\dots, y_k)$. Namely,
for any $G(k)$-cover $X \to Y$, there exist rational functions 
$\theta_1,\dots, \theta_k \in \CC(X)$ such that
\begin{enumerate}[(i)]
\item  the action of $G(k)$ on $\{\theta_1, \dots, \theta_k\}$ is given by
\begin{eqnarray*}
(\theta_1,\dots, \theta_k)^{\sigma} &=& (\theta_2, \theta_3,\dots, \theta_k, \theta_1) \\
(\theta_1,\dots, \theta_k)^{\tau_i} &=& (\theta_1,\dots,\theta_{i-1}, -\theta_i, - \theta_{i+1}, \theta_{i+2},\dots, \theta_k) \quad (i = 1,\dots, k-1)
\end{eqnarray*}
and
\item $\CC(X) = \CC(Y)(\theta_1,\dots, \theta_k)$, $\CC(Y) = \CC(X)^{G(k)}$.
\end{enumerate}
\end{rem}

\begin{prop}\label{prop:necGk} Let $\pi : X \to Y$ be a $G(k)$-cover. Let  $H(k)$ be the normal
subgroup of $G(k)$ generated by  $\tau_1,\dots, \tau_{k-1}$. Let $Z:=D_{H(k)}(X/Y)$ be the $\CC(X)^{H(k)}$-normalization of
$Y$. Then there exist $\varphi_1, \dots, \varphi_k \in \CC(Z)$ such that
\begin{enumerate}[\rm(i)]
\item $\varphi_i^{\sigma} = \varphi_{i+1},\ i = 1, \dots, k-1$, $\varphi_k^{\sigma}: = \varphi_1$ and
\item $\CC(X) = \CC(Z)(\sqrt {\varphi_1},\dots, \sqrt{\varphi_k})$.
\end{enumerate}
\end{prop}

\begin{proof} Let $\theta_i$ $(i = 1, \dots, k)$ be the rational functions in $\CC(X)$ as in Remark~\ref{rem-versal}. Put $\varphi_i = \theta_i^2$.
One can check that these functions satisfy the desired conditions.
\end{proof}

We now consider the construction problem for $G(k)$-covers. 

\begin{lem}\label{lem:funGk}
Let $f : Z \to Y$ be a $k$-cyclic cover of $Y$. Suppose that there exist rational functions
$\varphi_1, \dots, \varphi_k \in \CC(Z)$ such that

\begin{enumerate}[\rm(i)]
\item $\varphi_i^{\sigma} = \varphi_{i+1}, \quad \varphi^{\sigma}_k = \varphi_1$,
\item $\varphi_1\cdots\varphi_k \in \left (f^*\CC(Y)^\times \right )^2$, $\varphi_i \not\in \left (\CC(Z)^{\times}\right )^2$
and,
\item $\sqrt {\varphi_{i+1} }\not\in \CC(Z)(\sqrt{\varphi_1},\dots, \sqrt{\varphi_i})$.
\end{enumerate}

Then $K:= \CC(Z)(\sqrt{\varphi_1}, \dots, \sqrt{\varphi_k})$ is a $G(k)$-extension and the $K$-normalization
$X$  of $Y$ is a
$G(k)$-cover of $Y$ such that $D_{H(k)}(X/Y)= Z$.
\end{lem}

\begin{proof} By the uniqueness of the normalization and Galois theory, the statement is immediate.
\end{proof}

\begin{prop}\label{prop:sufGk}
Let $f : Z \to Y$ be a $k$-cyclic cover of $Y$. Suppose that $Z$ is smooth and there
exist distinct irreducible divisors $D_1, \dots, D_k$ on $Z$ such that
\begin{enumerate}[\rm(i)]
\item $D_i \sim D_j$ for any $i, j$, and 
\item $D_i^{\sigma} = D_{i+1}$ $(i = 1, \dots, k-1)$ and $D^{\sigma}_k = D_1$.
\end{enumerate}
Then there exists an $H(k)$-cover $g : X \to Z$ such that 
\begin{enumerate}[\rm(a)]
\item $f\circ g : X \to Y$ is a $G(k)$-cover with
$D_{H(k)}(X/Y) = Z$ and
\item $\Delta_g = \Supp(D_1+ \cdots + D_k)$.
\end{enumerate}
\end{prop}

\begin{proof} Let $\varphi_i$ be rational functions such that 
\[
(\varphi_i) = D_{i+1} - D_i \quad (i = 1, \dots, k-1), \qquad (\varphi_k) = D_1 - D_k.
\]
Then $\varphi_i^{\sigma} = \varphi_{i+1}, \quad \varphi^{\sigma}_k = \varphi_1$ and $\varphi_1\cdots \varphi_k$ is
a constant. 
Let us denote $K_i:=\CC(Z)(\sqrt{\varphi_1}, \dots, \sqrt{\varphi_i})$.

\begin{claim} 
Under the above conditions 
$\sqrt{\varphi_{i+1} }\not\in K_i$, $1 \le i \le k-2$.
\end{claim}

\medskip

\begin{proof}[Proof of Claim.] Let $X_i$ be the $K_i$-normalization of $Z$ and  let
$g_i : X_i \to Z$ be the covering morphism.  Since$(\varphi_i) = D_{i+1} - D_i$, 
we have $\Delta_{g_i} = 
\Supp(D_1+\cdots + D_{i+1})$. Suppose that $\sqrt {\varphi_{i_0}} \in K_{i_0}$ for 
some $i_0$. Then $X_{i_0} = X_{i_0-1}$.  Hence $\Delta_{g_{i_0}} = \Delta_{g_{i_0 + 1}}$, 
which contradicts the hypothesis on the $D_i$'s being different.
\end{proof}

By the Claim and Lemma~\ref{lem:funGk}, the $K_k$-normalization $X$ of $Z$
satisfies the conditions of the statement. 
\end{proof}


\section{Proof of Theorem \ref{thm:main1}}

We keep the notations as before. 

\begin{lem}\label{lem:Drat}
If $D$ is rational, there exists a $D_{2n}$-cover of $\Sigma_1$
branched at $2(\bar{L}_1 + \bar{L}_2) + n\bar{D}$ for any $n \ge 3$.
\end{lem}

\begin{proof}  By Lemma \ref{lem:preimD}, $f^*{\bar D} = D^+ + \sigma_f^*(D^+)$ and
$D^+ \sim \sigma_f^*(D^+)$. By Proposition \ref{prop:sufdi}, the assertion follows.
\end{proof}

\par\bigskip

\begin{lem}\label{lem:Dnonrat}{If $D$ is not rational, there exist no
$D_{2n}$-cover of $\Sigma_1$ branched at $2(\bar{L}_1 + \bar{L}_2) + n\bar{D}$
for any $n \ge 3$.
}
\end{lem}

\begin{proof}  By Lemma \ref{lem:preimD}, $f^*\bar{D}$ is irreducible if $D$ is not
rational. By Corollary \ref{cor:branch}, our assertion follows.
\end{proof}
\par\bigskip

Now Theorem \ref{thm:main1} is straightforward consequence of the following.

\begin{prop}\label{prop:main1}
Let $n\geq 3$.
\begin{enumerate}[\rm(i)]
\item\label{prop:main1i} If there exists a $D_{2n}$-cover, $n$ odd, of 
$\PP^2$ branched at $2(L_1 + L_2) + nD$, then there exists a
$D_{2n}$-cover of $\Sigma_1$ branched at 
$2(\bar{L}_1 + \bar{L}_2) + n\bar{D}$.
\item If there exists a $D_{2n}$-cover of $\Sigma_1$ branched at 
$2(\bar{L}_1 + \bar{L}_2) + n\bar{D}$, then there exists a $D_{2n}$-cover 
of $\PP^2$ branched at $2(L_1 + L_2) + nD$.
\end{enumerate}
\end{prop}

\begin{proof} \begin{enumerate}[\rm(i)]
\item Suppose that there exists a $D_{2n}$-cover branched at
$2(L_1 + L_2) + nD$.  The double cover $\beta_{2,\langle\tau\rangle}(\pi) : D(S/\PP^2) \to \PP^2$ is
branched at $2(L_1 + L_2)$. Hence the double cover $f : \Sigma_2 \to \Sigma_1$ is
the canonical resolution of $\beta_{1, \langle\tau\rangle}(\pi) : D(S/\PP^2) \to \PP^2$ (see \cite {horikawa},
for the canonical resolution of double covers). Let $\tilde {S}$ be the
$\CC(S)$-normalization of $\Sigma_1$. The surface $\tilde S$ is a $D_{2n}$-cover and we denote the covering morphism
by $\tilde {\pi} : \tilde {S} \to \Sigma_1$.  Then $\beta_{1, \langle\tau\rangle}(\tilde {\pi}) = f$.
Since $\Delta_{\pi} = L_1 + L_2 + D$,  we have $\Delta_{\tilde {\pi}} \subset
\bar {L}_1 + \bar {L}_2 + \bar {D} + E$. If $\tilde {\pi}$ is branched along $E$, 
then $\beta_{2, \langle\tau\rangle}(\tilde {\pi})$ is branched along $E$. But this this impossible
by Corollary~\ref{cor:branch}, since $f^*E$ is irreducible. Thus \eqref{prop:main1i} follows.

\item If there exists a $D_{2n}$-cover $\tilde {\pi} : \tilde S \to \Sigma_1$ branched 
at $2(\bar {L}_1 + \bar {L}_2) + n \bar D$, the Stein factorization of $\mu_P\circ \tilde {\pi}$ gives
the desired $D_{2n}$-cover.
\end{enumerate}
\end{proof}


\section{Proof of Theorem \ref{thm:main2}}
\label{sec-proof-ii}

Let $C := D+ L_1 + L_2$ be a curve of type II, let $q : \Sigma_1 \to \PP^2$ be the blowing-up at $P$ and let $\overline {D}, \overline{L}_1$ and $\overline{L}_2$ be the 
strict transforms introduced in \S\ref{prel}. Let $f_k : \Sigma_k \to \Sigma_1$ be a
$k$-cyclic cover branched at $k(\overline {L}_1 + \overline{L}_2)$.

\begin{lem}\label{lem:Gkexist} 
If $D$ is rational, then there exists a $G(k)$-cover $S$ of $\Sigma_1$ branched at
$k(\overline {L}_1 + \overline{L}_2) + 2\overline{D}$ with 
$D_{H(k)}(S/\Sigma_1) = \Sigma_k$.
\end{lem}

\begin{proof} By Lemma \ref{lem:preimD},  
$f_k^*\overline{D} = G_1+ \cdots + G_k$, $G_i \sim G_j$ and
$G_i^{\sigma} = G_{i+1}$ $(i =1, \dots, k-1)$, $G_k^{\sigma} = G_1$. Hence by Proposition \ref{prop:sufGk}, 
we have a $G(k)$ cover with the desired properties.
\end{proof}

\bigskip

\begin{lem}\label{lem:Gknonexist}
If $D$ is not rational, then there exists no $G(k)$-cover  
$p_k: S \to \Sigma_1$ of $\Sigma_1$ such that

\begin{enumerate}[\rm(i)]
\item\label{lem:Gknonexisti}  $p_k$ is branched at either 
$k(\overline {L}_1 + \overline {L}_2) + 2\overline {D}$ or
\item\label{lem:Gknonexistii} $D_{H(k)}(S/\Sigma_1) = \Sigma_k$ and $\beta_{1, H(k)} = f_k$.
\end{enumerate}
\end{lem}

\begin{proof} Suppose that  there exists a $G(k)$-cover  $p_k  : S \to \Sigma_1$ of 
$\Sigma_1$ satisfying the conditions \eqref{lem:Gknonexisti} and \eqref{lem:Gknonexistii}. 
Then by Proposition~\ref{prop:necGk}, there exist functions $\varphi_1, \dots, \varphi_k
\in \CC(D_{H(k)}(S/\Sigma_1))$ such that

\begin{enumerate}[(a)]
\item $\CC(S) = \CC(D_{H(k)}(S/\Sigma_1))(\sqrt{\varphi_1}, \dots, \sqrt{\varphi_k})$, and
\item $\varphi_i^{\sigma} = \varphi_{i+1} (i =1,\dots, k-1)$, 
$\varphi_k^{\sigma} = \varphi_1$. 
\end{enumerate} 

By Lemma \ref{lem:preimD}, we have $f^*\overline D = G_1+\cdots + G_\lambda$ 
$(\lambda | k, \lambda < k)$. Since $\Delta_{\beta_{2, H(k)}} = f^*\overline D$ or 
$f^*\overline D + \Delta_{0, k}$, the divisor of $\varphi_i$ is of the form
\[
(\varphi_i) = \sum_j c_{ij}G_j + c'_i\Delta_{0,k} + 2L_i,
\]
where $c_{ij}, c'_i \in \{ 0, 1\}$. Since  $\Delta_{0,k}$ is $\sigma$-invariant, 
$G_i^{\sigma} = G_{i+1}$ and $G_i^{\sigma^l}  = D_i$, we have 
\[
(\varphi_{i+l}) = (\varphi_i^{\sigma^l}) = \sum_j c_{ij}G_j + c'_i\Delta_{0,k} + 2L_i^{\sigma^l}.
\]
Hence $(\frac{\varphi_i}{\varphi_{i+l}}) = 2(L_i - L_i^{\sigma^l})$. As 
$D_{H(k)}(S/\Sigma_1) = \Sigma_k$ is simply connected, 
$\frac{\varphi_i}{\varphi_{i+l}} \in \left (\CC(D_{H(k)}(S/\Sigma_1))^{\times}\right )^2$. This implies that
$[\CC(S):\CC(D_{H(k)}(S/\Sigma_1))] \le 2^\lambda < 2^{k-1}$ which contradicts $\lambda<k$.
\end{proof}

\bigskip

\begin{cor}\label{cor:Gknonexist}
If $k$ is odd and $D$ is not rational, then there exists no $G(k)$-cover
$\pi : S \to \Sigma_1$ branched at $k(\overline{L}_1 + \overline{L}_2) + 2 \overline{D}$.
\end{cor}

\begin{proof} 
If $\beta_{1, H(k)}(\pi)$ is branched at $\overline{D}$, then the ramification index along $\overline{D}$ is an odd number $k_1$ $(k_1 | k, k_1 > 1)$.  Hence the ramification index
along $\overline D$ of $\pi$ is either $k_1$ or $2k_1$, either way different from 2.
Hence $\Delta_{\beta_{1, H(k)}(\pi)} = \overline {L}_1 + \overline {L}_2$.
\end{proof}

\begin{proof}[Proof of Theorem {\rm\ref{thm:main2}}]

Let $C$ be a curve of type II. 
\begin{itemize}
\item[\eqref{thm:main2i}] If $D$ is rational, then by Lemma \ref{lem:Gkexist}, there exists 
a $G(k)$-cover $\pi : S \to \Sigma_1$ of $\Sigma_1$. The Stein factorization of $q\circ\pi$
gives rise to the desired $G(k)$-cover.

\item[\eqref{thm:main2ii}] Suppose that $D$ is not rational and that there exists a
$G(k)$-cover branched at $k(\overline{L}_1 + \overline{L}_2) + 2D$. We denote it by 
$\overline{\pi} : \overline{S} \to \PP^2$. Let 
$q : \Sigma_1 \to \PP^2$ be a blowing-up at $P$ and let $S$  be the $\CC(\overline{S})$-normalization of $\Sigma_1$. 
$S$ is a $G(k)$-cover of $\Sigma_1$, and we denote its covering morphism by $\pi$.  By its definition, the branch locus
$\Delta_{\pi}$ of $\pi$ is either $\overline {L}_1 + \overline{L}_2 + \overline D$ or $\overline{L}_1 + \overline{L}_2 + 
\overline{D} + \Delta_{0,1}$.
 
\begin{lem}\label{lem:suplem1}{$\beta_{1, H(k)} : D_{H(k)}(S/\Sigma_1) \to \Sigma_1$ is branched along
$\overline{L}_1 + \overline{L}_2$
}
\end{lem}

\begin{proof} Since $S\setminus \pi^{-1}(\Delta_{0,1}) \to \Sigma_1$ is identified with $S\setminus \pi^{-1}(P) \to 
\PP^2\setminus \{P\}$, $\Delta_{\beta_1, H(k)}$ is either $\overline{L}_1 + \overline{L}_2$ or 
$\overline{L}_1 + \overline{L}_2 + \Delta_{0,1}$. If the latter case occurs, then the restriction of
$\Delta_{\beta_1, H(k)}$ to a general fiber of $\Sigma_1$ gives a $\ZZ/k\ZZ$-cover of $\PP^1$ branched at 
one point, which is impossible.

\end{proof}

\par\bigskip

By Lemma~\ref{lem:suplem1}, we have $D_{H(k)}(S/\Sigma_1) = \Sigma_k$ and $\beta_{1, H(k)}(\pi) = f_k$. Now by
Lemma~\ref{lem:Gknonexist}, we infer that $\pi : S \to \Sigma_1$ does not exist. Hence there exists no $G(k)$-cover
branched at $k(L_1 + L_2) + 2D$.
\end{itemize}
\end{proof}

\section{Fundamental groups of curves of type I and II}

Let $C$ be a rational curve of type I and degree $d$.
Let $P_i\in L_i\cap D$ ($P\neq P_i$) be either a singular point
of $D$ of type $\bbA_{2r_i}$ or a smooth point of $D$ (by convention of type $\bbA_0$,
i.e. $r_i=0$). We say that $C$ is nodal if the remaining singular points of $D$
(outside $P,P_1$ and $P_2$) are nodes. By the genus formula, the number of such
nodal points equals $d-2-r_1-r_2$.
We want to compute $\pi_1(\PP^2\setminus C)$. We follow the notation
of section \S\ref{prel}. Since $P\in C$, it is clear that this group is naturally
isomorphic to $G:=\pi_1(\Sigma_1\setminus (\overline{D}\cup\overline{L_1}\cup\overline{L_2}\cup\Delta_{0,1}))$.
Since
$$
f_{2|}:\Sigma_2\setminus (f_2^*\overline{D}\cup f_2^*\overline{L_1}\cup f_2^*\overline{L_2}\cup\Delta_{0,2})\to\Sigma_1\setminus (\overline{D}\cup\overline{L_1}\cup\overline{L_2}\cup\Delta_{0,1})
$$
is an unramified cyclic covering we have a short exact sequence:
\begin{equation}\label{suc-ex}
1\to
H:=\pi_1(\Sigma_2\setminus (f_2^*\overline{D}\cup f_2^*\overline{L_1}\cup f_2^*\overline{L_2}\cup\Delta_{0,2}))\to
G
\to\ZZ/2\ZZ\to 1.
\end{equation}
In order to compute $H$ we will use some facts about the geometry of $C$.
Note that $\Sigma_2\setminus (f_2^*\overline{L_2}\cup\Delta_{0,2})$
is isomorphic to $\CC^2$ and that the restriction of $p_2$ to 
$\Sigma_2\setminus (f_2^*\overline{L_2}\cup\Delta_{0,2})$, induces
a projection onto $\CC$. Thus, in order to compute $H$, we can apply braid monodromy
and the asymptotic version of Zariski-van Kampen theorem (see \cite{car:xx}).
The affine part of $\Sigma_2\setminus (f_2^*\overline{D}\cup f_2^*\overline{L_1}\cup f_2^*\overline{L_2}\cup\Delta_{0,2})$ is given by 
$\Sigma_2\setminus (f_2^*\overline{D}\cup f_2^*\overline{L_1})$ 
which is the union of a fiber and two sections of $p_2$.
The non-transversal fibers of $p_2$ with respect to $f_2^*\overline{D}$ are of three types:
\begin{enumerate}
\makeatletter\renewcommand{\theenumi}{T\arabic{enumi}}
\item\label{item-l1} $f_2^*\overline{L_1}$.
\item\label{item-l2} 
Vertical lines through double points of $f_2^*\overline{D}$ (which correspond to transversal
intersections of the irreducible components of $G_1,G_2$). This case does not happen
if $r_1+r_2=d-2$.
\item\label{item-l3} Vertical asymptotic lines to $f_2^*\overline{D}$ (which correspond to transversal intersections of the irreducible components of $G_1,G_2$ with $\Delta_{0,2}$).
There are $d-2$ such lines for each $G_i$.
\makeatletter\renewcommand{\theenumi}{\arabic{enumi}}
\end{enumerate}

Let us study the local behavior of braid monodromy around these lines:
\begin{itemize}
\item[\eqref{item-l1}] Each $G_i$
intersect at a point $P_1\in f_2^*\overline{L_1}$; their intersection number equals $2 r_1+1$.
Local braid monodromy equals the braid $\sigma_1^{2(2 r_1+1)}$, and the relations
are of type $\mu_i=\mu'_i$, where  $\mu_i,\mu'_i$ are some meridians
of $G_i$.

\item[\eqref{item-l2}] In this case, the local braid monodromy equals the braid 
$\sigma_1^2$, and the relations are of the type $[\mu_1,\mu_2]=1$, where  $\mu_i$ is
a meridian of $G_i$.

\item[\eqref{item-l3}] In this case, the local braid monodromy equals the braid 
$\sigma_1^{-2}$, and the relations are of the type $[\mu_1,\mu_2]=1$, where  $\mu_i$ is
a meridian of $G_i$.

\end{itemize}

\begin{prop}\label{prop-grupo} Let $C$ be a curve of type I. Then:

\begin{enumerate}[\rm(a)]
\item\label{prop-grupoa} If $C$ is rational and nodal, then $G$ is isomorphic to $K_2$, where
$$
K_2:=\!\!\langle \ell,x_1,x_2\mid
[x_1,x_2]=1, 
\ell^{-1} x_1 \ell= x_2,\ell^{-1} x_2 \ell= x_1
\rangle\!\!\cong\!\!
\langle y,x\mid\!\!
[x^2,y]= 
[y^2,x]=1
\rangle.
$$

\item\label{prop-grupoc} If $C$ is non-rational, nodal and it is a generic
element of a family $\{C_t\}_{t\in\Delta}$, such that $C_0$ is as in \eqref{prop-grupoa},
then its fundamental group is Abelian.
\end{enumerate}
\end{prop}

\begin{proof}\mbox{}
As for \eqref{prop-grupoa} we break down the proof into two parts:
\begin{enumerate}[(i)]
\item There is an epimorphism:
We apply a modified version of Zariski-van Kampen Theorem as in 
\S\ref{sec-example} to the curve $f_2^*\overline{L_1}\cup G_1\cup G_2$
(see Lemma~\ref{lem:preimD} for notations).
According to such method, the group $H$ is generated by three meridians 
$\ell_1,x_1,x_2$, where $\ell_1$ is a meridian of $f_2^*\overline{L_1}$, 
and $x_i$ is a meridian of $G_i$ (the latter contained in a \emph{generic} 
vertical line $F$ near $f_2^*\overline{L_1}$). The relations given by 
$f_2^*\overline{L_1}$ are of type $\ell_1^{-1} x_i \ell_1=w^{-1} x_i w$ for a 
certain word $w$ in $\ell_1,x_1,x_2$ (in fact it can be chosen in $x_1,x_2$).
The relations given by the lines of type \eqref{item-l2} and \eqref{item-l3} are 
commutators of a conjugate of $x_1$ and a conjugate of $x_2$. This means that 
there exists a surjection of $H$ onto the Abelian group generated by $\ell_1,x_1,x_2$. 
Using the exact sequence \eqref{suc-ex}, we can choose a meridian $\ell$ of $L_1$ such 
that $\ell^2=\ell_1$ and thus we are done.

\item It is an isomorphism:
The hypothesis implies the existence of a non-transversal fiber $F_1$ of type 
\eqref{item-l2} or \eqref{item-l3}. Note that the meridians $x_1,x_2$ in $F$ coincide
with the meridians $\mu_1,\mu_2$ described in \eqref{item-l2} and \eqref{item-l3}.
Since the braid group is $\BB_2\cong\ZZ$ we deduce $x_1$ and $x_2$ commute. 
Therefore $w$ commutes with $x_1,x_2$ and thus $\ell_1$ commutes with $x_1,x_2$.
\end{enumerate}

As for \eqref{prop-grupoc} we have an epimorphism of $K_2$ onto the fundamental group
of $C$ and it is easily seen that $x_1$ and $x_2$ are conjugated
in the subgroup generated by $x_1, x_2, \ell^2$, then $x_1=x_2$ and hence we are done.
\end{proof}

For curves of type II, we may proceed in a similar way. In this case one obtains
an exact sequence similar to \eqref{suc-ex}, where the right-most term is $\ZZ/k\ZZ$. 
Suitable replacements of $2$ by $k$ must be performed, e.g. $f_k^*(\overline{D})$
has now $k$ irreducible components $G_1,\dots,G_k$. For the
case \eqref{item-l3}, now the local braid is 
$(\sigma_1\dots\sigma_{k-2}\sigma_{k-1}^2\sigma_{k-2}\dots\sigma_1)^{-1}$ and
one meridian of one component $G_i$ commutes with meridians of the other components.
With the same ideas of Proposition~\ref{prop-grupo}, we can prove a similar
result.

\begin{prop}\label{prop-grupo2} Let $C$ be a curve of type II.

\begin{enumerate}[\rm(a)]
\item\label{prop-grupo2a} If $C$ is rational and nodal, then there exists
an epimorphism $G\twoheadrightarrow K_k$, where
$$
K_k:=\langle \ell,x \mid
[x,\ell^k]=[x,\ell^{-i}x\ell^i]=1, i=1,...,k-1
\rangle.
$$
\item\label{prop-grupo2b} If $C$ is rational, nodal, has real equations
and either  all the branches at $P$ are real or all nodes are real with real tangent lines, then the
epimorphism of \eqref{prop-grupo2a} is an isomorphism.

\item\label{prop-grupo2c} If $C$ is non-rational, nodal and it is a generic
element of a family $\{C_t\}_{t\in\Delta}$, such that $C_0$ is as in \eqref{prop-grupo2b},
then its fundamental group is $K_\lambda$ for some $\lambda|k$, $1\leq \lambda<k$.
\end{enumerate}
\end{prop}

\begin{proof}
The proof of \eqref{prop-grupo2a} follows along the lines of the epimorphism part in the proof of
Proposition~\ref{prop-grupo}\eqref{prop-grupoa}, changing $2$ by $k$ where necessary.
This way one obtains the following group:
\begin{equation}
\label{eq-kk}
\langle \ell,x_{1},\dots,x_{k} \mid
[x_{i},x_{j}]=1,\ell^{-1} x_{i} \ell= x_{{i+1}},
\ i,j=1,...,k\!\!\mod k
\rangle
\end{equation}

Finally, eliminating $x_i$, $i=2,...,k$ and setting $x=x_1$ provides the given 
presentation with two generators.

The arguments of \eqref{prop-grupo2b} differ from those in the proof of Proposition~\ref{prop-grupo}. The hypothesis implies that the $k$ sections are real and
intersect both each other and $\Delta_{0,k}$ transversally at real points with real 
branches (outside $f_k^*\overline{L_1}$ and $f_k^*\overline{L_2}$). Checking at every
non-transversal fiber of type \eqref{item-l2} or \eqref{item-l3}, one sees that the 
meridians involved in such relations are $x_1,\dots,x_k$ (see the proof of
Proposition~\ref{prop-grupo}\eqref{prop-grupoa}); in this case, we use for it that we can compute
the braid monodromy from a real picture. Since all the sections intersect each 
other, then $[x_i,x_j]=1$, $\forall i,j$. Therefore $w$ commutes with
$x_1,\dots,x_k$ and thus also $\ell_1$.

Let us prove \eqref{prop-grupo2c}. Note that, in this case, $G$ is a quotient of $K_k$.
Some commutation relations $[x_i,x_j]=1$ obtained in \eqref{prop-grupo2b} are replaced
by relations $x_i=x_j$. Each irreducible component $G_1,\dots,G_\lambda$ is associated
with one block of generators in $\{x_1,...,x_k\}$, and the $\frac{k}{\lambda}$ generators of each block are equal (they are meridians of the corresponding component). Let us denote these meridians by $y_1\dots,y_\lambda\in\{x_1,...,x_k\}$. In particular, $[y_i,y_j]=1$.
Since the $k$-cyclic Galois action permutes the irreducible components, and this action 
induces conjugation by $\ell$, we can order these meridians so that 
$\ell^{-1} y_i \ell= y_{i+1}$, $1\leq i<\lambda$, and $\ell^{-1} y_\lambda \ell= y_1$.
In other words $G\cong K_\lambda $.
\end{proof}


\end{document}